\documentclass[12pt,oneside]{article}

\usepackage[english]{babel}

\usepackage[letterpaper,top=2cm,bottom=2cm,left=3cm,right=3cm,marginparwidth=1.75cm]{geometry}

\usepackage{amsmath}
\usepackage{amsthm}
\usepackage{amssymb}
\usepackage{times}
\usepackage{tikz-cd}
\theoremstyle{plain}
\newtheorem{theorem}{Theorem}
\newtheorem{lemma}{Lemma}
\newtheorem{corollary}{Corollary}
\newtheorem{proposition}{Proposition}
\newtheorem{problem}{Problem}
\theoremstyle{definition}
\newtheorem{definition}{Definition}
\newtheorem{remark}{Remark}
\newtheorem{exam}{Example}

\usepackage{graphicx}
\usepackage[colorlinks=true, allcolors=blue]{hyperref}
\usepackage[english]{babel}
\usepackage{mathtools}
\usepackage{indentfirst}
\usepackage{xcolor}  
\usepackage{array}

\begin{document}

\title{An Upper Bound on the Number of Generalized Cospectral Mates of Oriented Graphs}
\author{\small Limeng Lin$^{{\rm a}}$$\quad\quad$ Wei Wang$^{\rm a}$\thanks{Corresponding author: wang\_weiw@xjtu.edu.cn}$\quad\quad$  Hao Zhang$^{{\rm b}}$
		\\
		{\footnotesize$^{\rm a}$School of Mathematics and Statistics, Xi'an Jiaotong University, Xi'an 710049, P. R. China}\\
	\footnotesize$^{\rm b}${ School of Mathematics, Hunan University, Changsha, 410082, P.R. China}		}
\date{}
	\maketitle

\begin{abstract}
This paper examines the spectral characterizations of oriented graphs. Let $\Sigma$ be an $n$-vertex oriented graph with skew-adjacency matrix $S$. Previous research mainly focused on self-converse oriented graphs, proposing arithmetic conditions for these graphs to be uniquely determined by their generalized skew-spectrum ($\mathrm{DGSS}$). However, self-converse graphs are extremely rare; this paper considers a more general class of oriented graphs $\mathcal{G}_{n}$ (not limited to self-converse graphs), consisting of all $n$-vertex oriented graphs $\Sigma$ such that $2^{-\left \lfloor \frac{n}{2}  \right \rfloor }\det W(\Sigma)$ is an odd and square-free integer, where $W(\Sigma)=[e,Se,\dots,S^{n-1}e]$ ($e$ is the all-one vector) is the skew-walk matrix of $\Sigma$. Given that $\Sigma$ is cospectral with its converse $\Sigma^{\rm T}$, there always exists a unique regular rational orthogonal $Q_0$ such that $Q_0^{\rm T}SQ_0=-S$. This study reveals that there exists a deep relationship between the level $\ell_0$ of $Q_0$ and the number of generalized cospectral mates of $\Sigma$.
More precisely, we show, among others, that the maximum number of generalized cospectral mates of $\Sigma\in\mathcal{G}_{n}$ is at most $2^{t}-1$, where $t$ is the number of prime factors of $\ell_0$.
Moreover, some numerical examples are also provided to demonstrate that the above upper bound is attainable. Finally, we also provide a criterion for the oriented graphs $\Sigma\in\mathcal{G}_{n}$ to be weakly determined by the generalized skew-spectrum ($\mathrm{WDGSS})$.

\end{abstract}

    \section{Introduction}
    Let $G$ be a simple graph on $n$ vertices with vertex set $V(G)=\{v_1,\ldots,v_n\}$ and edge set $E(G)$. The \emph{adjacency matrix} of $G$ is an $n$ by $n$ matrix $A(G)=\left [ a_{i,j} \right ]_{n\times n}$, where $a_{i,j} = 1$ if $v_i$ and $v_j$ are adjacent, and $a_{i,j} = 0$ otherwise. The \emph{spectrum} of $G$, denoted by ${\rm Spec}(G)$, consists of all the eigenvalues (including the multiplicities) of $A(G)$. The complement of the graph $G$, denoted by $\bar{G}$, is a graph with the same vertex set with $G$, while two vertices are adjacent in $\bar{G}$ if and only if they are non-adjacent in $G$. The \emph{generalized spectrum} of $G$ is the ordered pair $({\rm Spec}(G),{\rm Spec}(\bar{G}))$.

Two graphs $G$ and $H$ are \emph{cospectral} if they share the same spectrum, i.e., ${\rm Spec}(G)={\rm Spec}(H)$. A graph $G$ \textit{is determined by its spectrum} if any graph cospectral with $G$ is isomorphic to $G$, abbreviated as $\mathrm{DS}$. A graph $H$ is referred to as a \textit{cospectral mate} of $G$ if $H$ is cospectral with $G$ but non-isomorphic to $G$. In other words, a graph is $\mathrm{DGS}$ if it has no cospectral mate.

A long standing unsolved question in spectral graph theory is ``Which graphs are $\mathrm{DS}$?''. The problem originates more than sixty years ago from chemistry and was first raised by G\"{u}nthard and Primas~\cite{GP}. It is also closely related to many problems of central interests such as the graph isomorphism problem and the famous problem of Kac~\cite{KAC}: ``Can one hear the shape of a drum?". Generally speaking, it is very hard and challenging to show a graph to be DS and up to now, only very few families of graphs with special structures were shown to be DS. For more background and known results, we refer the reader to \cite{DH,DH1} and the references therein.

 Wang and Xu \cite{ref7} studied the above problem from the perspective of generalized spectrum. Two graphs $G$ and $H$ are \emph{generalized cospectral} if they share the same generalized spectrum, i.e., ${\rm Spec}(G)={\rm Spec}(H)$ and ${\rm Spec}(\bar{G})={\rm Spec}(\bar{H})$. A graph $G$ \textit{is determined by its generalized spectrum} if any graph generalized cospectral with $G$ are isomorphic to $G$, abbreviated as $\mathrm{DGS} $.
  Let $W(G):=\left [ e,Ae,\dots,A^{n-1}e\right ]$ be the \emph{walk matrix} of $G$, where $e$ is the all-one vector. Then it can be shown that $2^{-\left \lfloor \frac{n}{2}  \right \rfloor }\mathrm{det}W(G)$ is always an integer, see~ \cite{W3}. In~\cite{W4}, Wang proved the following
\begin{theorem}[Wang~\cite{W4}]\label{GDGS}
   If $2^{-\left \lfloor \frac{n}{2}  \right \rfloor }\det W(G)$ is odd and square-free, then $G$ is DGS.
\end{theorem}

The above spectral determination problem naturally extends to oriented graphs.
 Given a simple undirected graph $G$ with an orientation $\sigma$. An oriented graph, denoted by $\Sigma = (G,\sigma )$, is a directed graph obtained from $G$ by assigning a direction to each edge of $G$ according to $\sigma$. We call $G$ the \emph{underlying graph} of $\Sigma $. The set of all oriented graphs on $n$ vertices is denoted by $\mathfrak{D}_n$.

 The \emph{skew-adjacency matrix} of an oriented graph $\Sigma $, first introduced by Tutte~\cite{ref11}, is the $n\times n$ matrix $S({\Sigma})=(s_{ij})$, where

\begin{equation}
     s_{i,j}=\left\{\begin{matrix}
 1 &  \mathrm{if}\ (v_{i}, v_{j})\ \mathrm{is\ an\ arc};  \\
 -1 & \mathrm{if}\ (v_{j}, v_{i})\ \mathrm{is\ an\ arc};\\
  0& \mathrm{otherwise.}
\end{matrix}\right.
\end{equation}
It is easy to see that $S(\Sigma)$ is a skew-symmetric matrix. The \emph{skew-spectrum} of $\Sigma$ is defined as the spectrum of $S(\Sigma)$. Specifically, in our paper, we define that the generalized skew-spectrum of a graph $\Sigma$ as the spectrum of $S(\Sigma)$ and the spectrum of $J-I-S(\Sigma)$, where $J$ and $I$ denote the all-one matrix and identity matrix of $n$ order, respectively. The aforementioned notions of spectral determination of graphs can be extended to oriented graphs accordingly. An oriented graph ${\Sigma}$ is \emph{determined by the generalized skew-spectrum} (DGSS for short), if any oriented graph ${\Delta}$ having the same generalized skew spectrum as ${\Sigma}$ is isomorphic to ${\Sigma}$.


An oriented graph $\Sigma$ is \emph{self-converse} if it is isomorphic to its converse $\Sigma ^{\rm T}$, a graph obtained by reversing each directed edge in $\Sigma $.
Let \[W(\Sigma):=\left [ e,Se,\dots,S^{n-1}e\right ]\]
 be the \emph{skew-walk matrix} of $\Sigma$. We call an oriented graph $\Sigma$ \emph{controllable} if $\mathrm{det}W(\Sigma)\ne 0$.

 Qiu et al.~\cite{ref4} extended Theorem~\ref{GDGS} to self-converse oriented graphs.
\begin{theorem}[\cite{ref4}]\label{DGSS}
Let $\Sigma$ be a self-converse oriented graph of order $n$. If $2^{-\left \lfloor \frac{n}{2}  \right \rfloor }\mathrm{det}W(\Sigma)$
(which is always an integer) is odd and square-free, then $\Sigma$ is $\mathrm{DGSS}$.
\end{theorem}
Recently, Chao et al.~\cite{ref5} further strengthened the above theorem by relaxing the constraint that $2^{-\left \lfloor \frac{n}{2}  \right \rfloor }\mathrm{det}W(\Sigma)$ must be square-free.

In Qiu et al.~\cite{ref4} and Chao et al.~\cite{ref5}, the authors are mainly concerned with self-converse oriented graphs. However, Wissing~\cite{ref9} gives theoretical evidence which shows that almost no mixed graphs are self-converse, i.e., the fraction of self-converse mixed graphs tends to zero as the order $n$ of graphs goes to infinity, indicating that self-converse graphs are rare. Therefore, in this paper, we shall consider a more general family of oriented graphs (not limited to self-converse ones), defined as follows:
\begin{center}
    $\mathcal{G}_{n}=\left \{ \Sigma\in {\mathfrak{D}_n}:\, 2^{-\left \lfloor \frac{n}{2}  \right \rfloor }\det W(\Sigma)~\text{is an odd and square-free integer}\right \}.$

\end{center}

Note that $\Sigma$ and $\Sigma^{\rm T}$ always have the same generalized skew-spectrum. It follows that a non-self-converse oriented graph $\Sigma$ always has a generalized cospectral mate $\Sigma^{\rm T}$, implying that every non-self-converse oriented graph is not $\mathrm{DGSS}$. Thus, in order for an oriented graph $\Sigma$ to be DGSS, it is necessary that it is self-converse. Due to this reason,  Qiu et al.~\cite{ref4} introduced the following notion. An oriented graph is called \textit{weakly determined by the generalized skew spectrum } $(\mathrm{WDGSS})$ if any oriented graph having the same generalized skew spectrum as $\Sigma$ is either isomorphic to $\Sigma$ or $\Sigma^{\rm T}$. It is natural to consider the following
 \begin{problem} [Qiu et al.~\cite{ref4}]\label{P1}
      Which oriented graphs are $\mathrm {WDGSS}$?
 \end{problem}
 In this paper, we provide a criterion for an oriented graph $\Sigma\in\mathcal{G}_{n}$ to be $\mathrm {WDGSS}$.
Moreover, for a non-$\mathrm {WDGSS}$ oriented graph $\Sigma$, we provide an upper bound for the number of its generalized cospectral mates, which is attainable as illustraed by some numerical examples.

To state our main results, we need some notions and notations.
 Given that $\Sigma$ and $\Sigma^{\rm T}$ are generalized cospectral, it follows that there exists a unique regular rational orthogonal matrix $Q_{0}$ satisfying
 \begin{equation}\label{EE1}
 Q_{0}^{\rm T}S(\Sigma)Q_{0}=S(\Sigma^{\rm T}),
 \end{equation}
 provided $\Sigma$ is controllable; see Lemma~\ref{QSQ}. Henceforth, we use $\ell_{0}$ to denote the level of matrix $Q_{0}$ (see Definition~\ref{level}).
 
 Let $\Sigma\in \mathcal{G}_{n}$ and $Q_{0}$ be defined as Eq.~\eqref{EE1} with level $\ell_{0}$. Then $\ell_0$ is the product of some distinct odd primes (see Remark~\ref{odd1}), i.e., $\ell_0=p_1p_2\cdots p_t$, where $p_i$'s are distinct odd primes for some $t\geq 0$. 
 
 The main result of the paper is the following theorem, which shows that the number of generalized cospectral mates of an oriented graph $\Sigma\in \mathcal{G}_{n}$ can be upper bounded by the number of factors of $\ell_0$ minus one.
 \begin{theorem}\label{2t-1}
  Let $\Sigma\in \mathcal{G}_{n}$. Suppose that $\ell_0=p_1p_2\cdots p_t$, where $p_i$'s are distinct odd primes for some $t\geq 0$. Then number of the generalized cospectral mates of $\Sigma$ is at most $2^{t}-1$.
 \end{theorem}
 \begin{remark}
   The upper bound in Theorem~\ref{2t-1} is attainable, see Examples 1-3 in Section~\ref{sec6}.
 \end{remark}

 \begin{remark}
   By Theorem~\ref{2t-1}, if $t=0$ (or equivalently $\ell_0=1$), then $\Sigma$ is $\mathrm{DGSS}$, since it has no cospectral mate; if $t=1$ (or equivalently $\ell_0=p_1$ is an odd prime), then the unique cospectral mate of $\Sigma$ is its converse $\Sigma^{\rm T}$, which shows that $\Sigma$ is $\mathrm {WDGSS}$.
 \end{remark}

 The proof of Theorem~\ref{2t-1} is based on the following Theorems \ref{l mid l0} and \ref{Q=QP}, which are of independent interest. Theorem \ref{l mid l0} reveals a direct relation between the level $\ell_{0}$ of $Q_{0}$ and the level $\ell$ of $Q\in\mathcal{Q}(\Sigma)$, where $\mathcal{Q}(\Sigma)$ denotes the set of all regular rational orthogonal matrices $Q$ such that $Q^{\rm T}S(\Sigma)Q$ is an oriented graph.
\begin{theorem}\label{l mid l0}
    Let $\Sigma\in \mathcal{G}_{n}$ and $Q_{0}$ be defined as in Eq.~\eqref{EE1} with level $\ell_{0}$. Suppose the level of $Q\in \mathcal{Q}(\Sigma)$ is $\ell$. Then $\ell\mid \ell_{0}$.
\end{theorem}
\begin{remark} It is obvious to see when $\Sigma$ is self-converse, then $\ell_0=1$ and every $Q\in \mathcal{Q}(\Sigma)$ has level one and hence is a permutation matrix. Therefore, every self-converse oriented $\Sigma\in\mathcal{G}_{n}$ is $\mathrm{DGSS}$. That is, Theorem~\ref{l mid l0} includes Theorem~\ref{DGSS} as a special case.
\end{remark}

The following theorem shows that two oriented graphs obtained from a controllable oriented graph via two regular rational orthogonal matrices of the same level, must be isomorphic.

\begin{theorem}\label{Q=QP}
     Let $\Sigma, \ \Delta,\ \Gamma \in \mathcal{G}_{n}$ be pairwise generalized cospectral. Suppose that there exist two regular rational orthogonal matrices $Q$ and $ Q_{1}\in\mathcal{Q}(\Sigma)$ with the same level $\ell$ such that
    $Q^{\rm T}S(\Sigma)Q=S(\Delta),$ and
    $Q_{1}^{\rm T}S(\Sigma)Q_{1}=S( \Gamma )$. Then $ \Delta$ is isomorphic to $\Gamma $.
\end{theorem}

Theorem~\ref{2t-1} implies that the count of cospectral mates is inherently constrained by the value of $\ell_{0}$. Moreover, for graphs $\Sigma\in \mathcal{G}_{n}$ with $\ell_{0}$ being the product of at most three prime factors, we can efficiently identify all cospectral mates of $\Sigma$ and elucidate the relationship among these cospectral graphs, based on an algorithm proposed in \cite{ref3}.

The rest of the paper is organized as follows. In Section~\ref{sec2}, we give some preliminary results that will be needed later in the paper. In Section~\ref{sec3}, we give the proof of Theorem~\ref{l mid l0}. In Section~\ref{sec4}, we present the proofs of Theorems \ref{2t-1} and \ref{Q=QP}. In Section~\ref{sec5}, we give a criterion for determining whether an oriented graph $\Sigma\in \mathcal{G}_{n}$ is $\mathrm{WDGSS}$. In Section~\ref{sec6}, we give some examples to illustrate Theorem \ref{2t-1}. Conclusions are given in Section~\ref{sec7}.

    \section{Preliminaries}\label{sec2}
    In this section, we shall give some preliminary results that will be needed later in the paper. If there is no confusion, we simply write $S = S(\Sigma)$ and $W = W(\Sigma)$.
\subsection{Main strategy}
    In this subsection, we shall describe the main strategy for showing a graph to be $\mathrm{DGSS}$. Recall a \emph{rational orthogonal matrix} $Q$ is an orthogonal matrix with rational entries, and it is called regular if $Qe=e$.

     The following lemma shows that a regular rational orthogonal matrix $Q$ can be employed to establish a relationship between two generalized cospectral oriented graphs, and whenever $\Sigma$ is controllable, the matrix $Q$ is unique. It serves as a pivotal element throughout this paper.
\begin{lemma}[\cite{ref4},\cite{ref7}] \label{QSQ}
    Let $\Sigma$ and $\Delta$ be two oriented graphs of order $n$. Suppose that $\Sigma$ is controllable. Then $\Sigma$ and $\Delta$
have the same generalized skew spectrum if and only if there exists a unique rational orthogonal matrix $Q$ such that
\begin{equation}
  Q^{\rm T}S(\Sigma)Q=S(\Delta),
\end{equation}
where $Q$ is regular, that is $Qe=e$.
\end{lemma}
Define
\begin{center}
$\mathcal{Q}(\Sigma) :=\cup_{\Delta\in \Delta_{n}} \left \{Q
\in \mathrm{RO}_n(\mathbb{Q} ): Q^{\rm T}S(\Sigma)Q = S(\Delta ) \right \} $
\end{center}
where $\mathrm{RO}_n(\mathbb{Q})$ denotes the set of all orthogonal regular matrices with rational entries of order $n$ and $\Delta_{n}$ denotes the set of all oriented graphs cospectral with $\Sigma$.
\begin{definition}[\cite{ref7}]\label{level}
     {\rm The level of a regular rational orthogonal matrix $Q$, denoted by $\ell(Q)$ or $\ell$, is the smallest
positive integer $k$ such that $kQ$ is an integral matrix.}
\end{definition}
\begin{lemma}[\cite{ref4}]
    Let $\Sigma$ be an oriented graph with $\mathrm{det}W(\Sigma)\ne 0$.
Then $\Sigma$ is determined by its generalized skew spectrum if and only if $\mathcal{Q}(\Sigma)$ contains only permutation matrices.
\end{lemma}
 It is obvious that a rational orthogonal matrix $Q$ with $Qe = e$ is a permutation matrix if and
only if $\ell(Q)=1$. Therefore, the above lemma can be sated in an alternative way.
\begin{corollary}[\cite{ref4}]\label{l=1}
    Let $\Sigma$ be an oriented graph with $\det W(\Sigma)\ne 0$. Then $\Sigma$ is $\mathrm{DGSS}$ if and only if the level $\ell$ of every $Q\in\mathcal{Q}(\Sigma)$ is 1.
\end{corollary}
The regular rational orthogonal matrix $Q_{0}$ such that $Q_{0}^{\rm T}S(\Sigma)Q_{0}=S(\Sigma^{\rm T})$ plays a key role in this paper. Of course, $Q_{0}$ must belong to $\mathcal{Q}(\Sigma)$. We record this in the following
\begin{proposition}
    Let $\Sigma$ be a controllable oriented graph. Then there exists a unique regular $Q_{0}\in {\mathcal{Q}(\Sigma)}$ with level $\ell_{0}$ such that
    \eqref{EE1} holds. Moreover, if $\Sigma$ is non-self-converse, then $Q_{0}$ is not a permutation matrix.
\end{proposition}

\subsection{Smith normal form}
An $n \times  n$ matrix $U$ with integer entries is called \emph{unimodular} if $\det U = \pm 1$. For every integral matrix $M$ with full rank, there exist two unimodular matrices $U$ and $V$ such that $M=USV$, where the matrix $S={\rm diag}(d_{1}, d_{2}, \dots, d_{n})$ is known as the \emph{Smith normal form (\rm{SNF} for short)} of $M$ and $d_{i}$ is called as the \emph{ $i-$th invariant factors} of $M$.
 The next lemma shows the relationship between the level $\ell$ and the last invariant factors of $W(\Sigma)$.
\begin{lemma}[\cite{ref7}]\label{l dn}
    Let $\Sigma$ be an oriented graph. Let $Q\in \mathcal{Q}(\Sigma)$ with level $\ell$. Then $\ell \mid d_{n}$, where $d_{n}$ is the $n$-th invariant factor of $W(\Sigma)$.
\end{lemma}

The following lemma gives the SNF of the skew-walk matrix $W(\Sigma)$ for a $\Sigma\in{\mathcal{G}_{n}}$.
\begin{lemma}[\cite{ref7}]\label{LQ}
        Let $\Sigma\in{\mathcal{G}_{n}}$. Then the SNF of $W(\Sigma)$ is as follow:

     \[{\rm diag}( \underbrace{1,1,\dots,1}_{\left \lceil \frac{n}{2}\right \rceil},\underbrace{2,2,\dots,2,2b}_{\left \lfloor \frac{n}{2}  \right \rfloor} ),\]

  where $b$ is an odd and square-free integer.
    \end{lemma}

\subsection{The level of $Q\in \mathcal{Q}(\Sigma)$ is odd}\label{PQ}
The following theorem will be needed later in the paper, the original proof of which is quite involved in~\cite{ref4}. Here, we provide a simplified alternative proof, based on the ideas from \cite{ref6}.
\begin{theorem}[\cite{ref4}]\label{lnot2}
     Let $\Sigma\in \mathcal{G}_{n}$ and $Q\in\mathcal{Q}(\Sigma)$ with level $\ell$. Then $2\nmid \ell$.
\end{theorem}
Before proving the above theorem, we need the following lemmas.
\begin{lemma}[cf.~\cite{W4}]
    Let $A$ be an $n\times n$ integral skew-symmetric matrix. If $A^{2} \equiv \mathbf{0}   \pmod{2}$, then $Ae \equiv 0\pmod{2}$.
\end{lemma}
\begin{proof}
    Let $A = (a_{i,j})$. Then the $(i, i)$-th entry of $A^{2}$ is
    \begin{center}
        $(A^{2})_{i,i}=\sum_{s=1}^{n}a_{i,s}a_{s,i}=-\sum_{s=1}^{n}a_{i,s}^{2}\equiv \sum_{s=1}^{n}a_{i,s}^{2}\equiv \sum_{s=1}^{n}a_{i,s}\pmod{2}$,
    \end{center}
which shows that $Ae \equiv 0\pmod{2}$.
\end{proof}
\begin{lemma}[cf. \cite{W4}]
     Let $\Sigma\in \mathcal{G}_{n}$ . Let $\varphi(x) = x^{n} + c_{1}x^{n-1} +\dots+ c_{n-1}x + c_{n}$ be the characteristic polynomial of graph $\Sigma $. Denote
     \begin{center}
         $M=M(\Sigma ):=\left\{\begin{matrix}
 S^{\frac{n}{2} }+c_{2}S^{\frac{n-2}{2}}+\dots+c_{n-2}S+c_{n}I,\mathit{if\ n\ is\ even;} \\
S^{\frac{n+1}{2} }+c_{2}S^{\frac{n-1}{2}}+\dots+c_{n-3}S^{2}+c_{n-1}S,\mathit{if\ n\ is\ odd.}
\end{matrix}\right.$
     \end{center}
     Then $Me\equiv 0 \pmod{2}$.
\end{lemma}
\begin{proof}
  We only prove the case that $n$ is even, the case that $n$ is odd can be proved in a
similar way.
For $\Sigma \in \mathcal{G}_{n}$, we know that $c_{i}=0$ when $i$ is odd by Sach's Coefficients Theorem. By Cayley-Hamilton Theorem, we get $\varphi(S) = S^{n} + c_{2}S^{n-2} +\dots+ c_{n-2}S^{2} + c_{n}I=0$. Thus, we can obtain
\begin{align*}
   M^{2}&=(S^{\frac{n}{2} }+c_{2}S^{\frac{n-2}{2}}+\dots+c_{n-2}S+c_{n}I)^{2}\\
  &\equiv S^{n} + c_{2}^{2}S^{n-2} +\dots+ c_{n-2}^{2}S^{2} + c_{n}^{2}I\\
  &\equiv S^{n} + c_{2}S^{n-2} +\dots+ c_{n-2}S^{2} + c_{n}I\\
   &\equiv0\pmod{2}.
\end{align*}
It follows that $Me\equiv 0 \pmod{2}$. The proof is complete.
\end{proof}
\begin{proof}[Proof of Theorem~\ref{lnot2}]
    Let $\Delta$ be any graph cospectral with $\Sigma$, and let $S_{1}$ represent the skew-adjacency matrix of $\Delta$. Then they have the same characteristic polynomial. Define $M(\Delta)$ as follows:
    \begin{center}
        $M(\Delta):=\left\{\begin{matrix}
 S_{1}^{\frac{n}{2} }+c_{2}S_{1}^{\frac{n-2}{2}}+\dots+c_{n-2}S_{1}+c_{n}I,\mathit{if\ n\ is\ even;} \\
S_{1}^{\frac{n+1}{2} }+c_{2}S_{1}^{\frac{n-1}{2}}+\dots+c_{n-3}S_{1}^{2}+c_{n-1}S_{1},\mathit{if\ n\ is\ odd.}
\end{matrix}\right.$
    \end{center}
 Hence, we have $M(\Delta)e\equiv0\pmod{2}$.
Let $k=\left \lceil \frac{n}{2}  \right \rceil $. We define two new matrices
\begin{equation}
  \bar{W}(\Sigma)=\left [e,Se,S^{2}e,\dots,S^{k-1}e,\frac{M(\Sigma)e}{2}, \frac{SM(\Sigma)e}{2},\dots,\frac{S^{n-r-1}M(\Sigma)e}{2} \right ] ,
\end{equation}
and
\begin{equation}
  \bar{W}(\Delta)=\left [e,S_{1}e,S_{1}^{2}e,\dots,S_{1}^{k-1}e,\frac{M(\Delta)e}{2}, \frac{S_{1}M(\Delta)e}{2},\dots,\frac{S_{1}^{n-r-1}M(\Delta)e}{2} \right ],
\end{equation}
 which are obtained from the matrices $W(\Sigma )$ and $W(\Delta)$, respectively, by some modifications.
 It is easy to know $\mathrm{det}\bar{W}(\Sigma)=2^{-\left \lfloor \frac{n}{2}  \right \rfloor}\mathrm{det}W(\Sigma) $. By $Q^{\rm T}SQ=S_{1}$ and $Q^{\rm T}e=e$, we have $Q^{\rm T}M(\Sigma )=M(\Delta)$ and hence $Q^{\rm T}\bar{W}(\Sigma)=\bar{W}(\Delta)$.  It follows from Lemma \ref{l dn} that $\ell \mid \det(\bar{W}(\Sigma)) $. By the definition of $\Sigma\in \mathcal{G}_{n}$, we have $2^{\left \lfloor \frac{n}{2}  \right \rfloor} \parallel \det W(\Sigma)$, which means $2$ does not divide $\det\bar{W}(\Sigma)$. ($a^{k}\parallel b$ means $a^{k}\mid b$ but $a^{k+1}\nmid b$.) Thus, it follows $2\nmid \ell$. This proves the theorem.
\end{proof}
\begin{remark}\label{odd1}
    Note that $Q_0\in \mathcal{Q}(\Sigma)$. It follows that $\ell_0\mid d_{n}$ and $2\nmid \ell_0$. Hence, $\ell_0$ is the product of some distinct odd primes whenever $\ell_0\neq 1$.
\end{remark}

\section{Proof of Theorem~\ref{l mid l0} }\label{sec3}
In this section, we present the proof of Theorem \ref{l mid l0}. In the subsequent discussions, we focus on the case that $p$ is an odd primes. The approach employed for the proof of Theorem~\ref{l mid l0} was inspired by the method presented by Chao et al.~\cite{ref5}, in which the proof is valid for self-converse graphs.
 For an integral matrix $M$, we use ${\rm rank }\,M$ and ${\rm rank }_p\,M$ to denote respectively the rank of $M$ over $\mathbb{Q}$ and over the finite field $\mathbb{F}_p$. Write $\bar{Q}_{0}=\ell_{0}Q_{0}$ in what follows.

We need several lemmas below.

\begin{lemma}[\cite{ref4}]
    Suppose that $\mathrm{rank}_{p} W(\Sigma)=r$. Then the first $r$ columns of $W(\Sigma)$ consist of a basis of the column space of $W(\Sigma)$, over $\mathbb{F}_p$
\end{lemma}

\begin{lemma}\label{z}
Let $\Sigma\in \mathcal{G}_{n}$ and $Q\in \mathcal{Q}(\Sigma)$ with level $\ell$. Suppose that $p$ is a prime divisor of $\ell$. Then there exists an integral vector $z\not\equiv 0 \pmod{p}$ such that
$W^{\rm T}z\equiv 0\pmod{p}$.
\end{lemma}
\begin{proof}

Let $\bar{Q} = \ell Q$. It follows from $Q^{\rm T}S(\Sigma)Q=S(\Delta)$ and
$Q^{\rm T}e=e$ that $Q^{\rm T}W(\Sigma) = W(\Delta)$, and hence $W(\Sigma)^{\rm T}\bar{Q}=
\ell W(\Delta)$. The definition of $\ell$ implies that there must exist some column $z$ of $\bar{Q}$ that is nonzero over $\mathbb{F}_{p}$. As $p\mid \ell$, we get $W(\Sigma)^{\rm T}z \equiv 0 $ and $z^{\rm T}z = \ell^{2}\equiv 0 \pmod{p^{2}}$. This lemma is proven.
\end{proof}
\begin{lemma}[\cite{ref4}]\label{sv=cv}
     Under the conditions of Lemma \ref{z}, then $S^{T}z\equiv\lambda_{0}z\pmod{p}$, for some integer $\lambda_{0}$.
\end{lemma}
\begin{lemma}\label{Sv=0}
    Let $\Sigma \in \mathcal{G}_{n}$ and $Q_{0}$ be defined as in \eqref{EE1} with level $\ell_{0}$. Suppose that $p$ is not a prime divisor of $\ell_{0}$. If $W^{T}z=0$, then $S^{T}z=0$ over $\mathbb{F}_{p}$.
\end{lemma}
\begin{proof}
By Lemma \ref{sv=cv}, we have $S^{\rm T}z \equiv \lambda_{0}z  \pmod{p}$ for some integer $\lambda_{0}$, where $z\not\equiv0\pmod{p}$. Next, we claim $\lambda_{0} \equiv 0  \pmod{p}$.

Note that ${Q}_{0}^{\rm T}SQ_{0}=S^{\rm T}$. Then we have
      \begin{equation}\label{SQz}
          S(\bar{Q}_{0}z)= \bar{Q}_{0}S^{\rm T}z\equiv \lambda_{0}(\bar{Q}_{0}z) \pmod{p}.
      \end{equation}
     Furthermore,
      \begin{align*}
          e^{\rm T}(S^{\rm T})^{k}(\bar{Q}_{0}z)&=(-1)^{k}e^{\rm T}S^{k}(\bar{Q}_{0}z)\\&\equiv (-1)^{k}\lambda_{0}^{k}e^{\rm T}(\bar{Q}_{0}z)\\&=(-1)^{k}\lambda_{0}^{k}\ell_{0}e^{\rm T}z\equiv0\pmod{p},
      \end{align*}
       for integers $1\le k\le n-1 $, where the reason for the last equality is  $Q_{0}e=e$ and $W^{\rm T}z\equiv0\pmod{p}$. It follows $W^{\rm T}(\bar{Q}_{0}z)\equiv0\pmod{p}$, and hence we get $z$ and $\bar{Q}_{0}z$ are linearly dependent over $\mathbb{F}_{p}$, as $\mathrm{rank}_{p} W=n-1$. Then there exists some integer $\lambda_{1}\not\equiv0\pmod{p}$ such that $\bar{Q}_{0}z \equiv \lambda_{1}z  \pmod{p}$. In fact, since $p\nmid \ell_{0}$ , then $\mathrm{rank}_{p}\bar{Q}_{0}=n$. Suppose $\lambda_{1}\equiv0\pmod{p}$, it follows that the nullspace of $\bar{Q}_{0}$ only have zero solution over $\mathbb{F}_{p}$, which implies $z\equiv0\pmod{p}$; a contradiction.

       Combining  $Sz\equiv-\lambda_{0}z\pmod{p}$ and  $\bar{Q}_{0}z \equiv \lambda_{1}z  \pmod{p}$ with Eq.~(\ref{SQz}), then we have
       \begin{center}
           $\lambda_{0}(\lambda_{1}z)\equiv\lambda_{0}(\bar{Q}_{0}z)\equiv S(\bar{Q}_{0}z)\equiv S\lambda_{1}z\equiv-\lambda_{1}\lambda_{0}z \pmod{p}.
           $
       \end{center}
       Therefore, we get $\lambda_{0}\equiv0\pmod{p}$ as $p$ is odd.
       This completes the proof.
\end{proof}
\begin{lemma}\label{key}
 Let $\Sigma$ be an oriented graph and $Q_{0}$ be defined as in \eqref{EE1} with level $\ell_{0}$.
 Then $\bar{Q}_{0}S^{k}e= (-1)^{k}\ell_{0}S^{k}e$ for any $k \ge  0$.
\end{lemma}
\begin{proof}
    Since $Q_{0}^{\rm T}SQ_{0}=S^{\rm T}$ and $Q_{0}e=e$, it is easy to get  $\bar{Q}_{0}S^{k}e= (-1)^{k}\ell_{0}S^{k}e$.
\end{proof}
\begin{lemma}[\cite{ref5}]
    Let $Q_{0}$ be a regular orthogonal matrix such that $Q_{0}^{\rm T}SQ_{0} = S^{\rm T}$. Suppose that $W$ is non-singular. Then $Q_{0}$ is symmetric.
\end{lemma}
\begin{lemma}\label{rankQ0}
   Under the conditions of Lemma \ref{Sv=0}, then
   \begin{center}
        $\mathrm{rank}(\ell_{0}I-\bar{Q}_{0})=\mathrm{rank}(\ell_{0}I+\bar{Q}_{0})=\frac{n}{2} $ for even $n$,
   \end{center}
   and
   \begin{center}
       $\mathrm{rank}(\ell_{0}I-\bar{Q}_{0})=\frac{n-1}{2}$and $ \mathrm{rank}(\ell_{0}I+\bar{Q}_{0})=\frac{n+1}{2} $ for odd $n$,
   \end{center}
   over $\mathbb{F}_{p}$.
\end{lemma}
\begin{proof}
    First, consider the case when $n$ is even.
    By Lemma \ref{key}, for odd $i$, we have $(\ell_{0}I+\bar{Q}_{0})S^{i}e=0$, i.e., $S^{i}e\in\mathrm{ker}(\ell_{0}I+\bar{Q}_{0})$, since $Se,S^{3}e,\dots,S^{n-1}e$ are all linearly independent, then we get $\mathrm{rank}(\ell_{0}I+\bar{Q}_{0})\le n/2 $. Similarly, for even $i$, $\mathrm{rank}(\ell_{0}I-\bar{Q}_{0})\le n/2 $. It follows
    \begin{align*}
        n&=\mathrm{rank}_{p}(\ell_{0}I)\le\mathrm{rank}_{p}(\ell_{0}I+\bar{Q}_{0})+\mathrm{rank}_{p}(\ell_{0}I-\bar{Q}_{0})\\
        &\le\mathrm{rank}(\ell_{0}I+\bar{Q}_{0})+\mathrm{rank}(\ell_{0}I-\bar{Q}_{0})\le\frac{n}{2}+\frac{n}{2}=n.
    \end{align*}
   Therefore the first assertion is obtained. The case that $n$ is odd can be proved in a similar way.
\end{proof}
\begin{remark}
    \textup{Let $V_{\lambda}$ denote the eigenspace of $\bar{Q}_{0}$ corresponding to the eigenvalue $\lambda$ over $\mathbb{F}_{p}^{n}$. By the above Lemma, we have dim$V_{-\ell_{0}}=\frac{n}{2}$ and dim$V_{\ell_{0}}=\frac{n}{2}$ when $n$ is even; we have dim$V_{-\ell_{0}}=\frac{n-1}{2}$ and dim$V_{\ell_{0}}=\frac{n+1}{2}$ when $n$ is odd. }
\end{remark}
\begin{lemma}\label{vv=0}  Under the conditions of Lemma \ref{Sv=0}.  If $W^{T}z=0$, then $z^{\rm T}z\ne 0$ over $\mathbb{F}_{p}$.
\end{lemma}
\begin{proof}
   We restrict our proof to the case where $n$ is even since the case where $n$ is odd can be demonstrated analogously. Let $V$  denote $\mathbb{F}_{p}^{n}$. Firstly, we define a space as follows:
    \begin{center}
        $V_{0}:=Col_{\mathbb{F}_{p}}(W)=\mathrm{span}\left \langle e,Se,\dots,S^{n-2}e \right \rangle \subset V=\mathbb{F}_{p}^{n}. $
    \end{center}
     Note that $\mathrm{rank}_{p}(W)=\mathrm{rank}_{p}(W^{\rm T})=n-1$, we get dim$V_{0}^{\perp}=1$, and hence $V_{0}^{\perp}=\mathrm{span}\left \langle z \right \rangle $.  By the Lemma \ref{rankQ0}, it follows that $e,S^{2}e,\dots,S^{n-2}e$ are the $\frac{n}{2} $ eigenvectors of $\bar{Q}_{0}$ corresponding to $\ell_{0}$, whereas $Se,S^{3}e,\dots,S^{n-3}e$ are the $\frac{n}{2}-1$ eigenvetors of $\bar{Q}_{0}$ corresponding to $-\ell_{0}$.  Note that dim$V_{-\ell_{0}}=\frac{n}{2}$. Thus, there exists an $\alpha \in V\setminus V_{0}$ such that $\bar{Q}_{0}\alpha=-\ell_{0}\alpha$.
     It follows $\mathbb{F}_{p}^{n}=V=V_{0}\oplus\left \langle\alpha  \right \rangle$.

    On the contrary, suppose that $z^{\rm T}z=0$. Then $z^{\rm T}\left [ W,z \right ] =0$. Note that $\mathrm{rank}_{p}(W)=n-1$. It is not difficult to see $\mathrm{rank}_{p}\left [ W,z \right ]=n-1$, since otherwise we obtain $z=0$. It implies that $z\in V_{0}$ and let
        $z= {\textstyle \sum_{i=0}^{n-2}} a_{i}S^{i}e$,
    where $a_{n-2}\ne 0$.
    In fact,  suppose $z= {\textstyle \sum_{i=0}^{n-3}} a_{i}S^{i}e$. It follows from Lemma \ref{Sv=0} that $Sz={\textstyle \sum_{i=0}^{n-3}} a_{i}S^{i+1}e=0$. Since $Se,\dots,S^{n-2}e$ are linearly independent, then $a_{i}=0$ for $1\le i\le n-2$ and hence $z=0$; a contradiction. Therefore $a_{n-2}\ne 0$.

     Recall Lemma \ref{Sv=0}, we know $\bar{Q}_{0}z=\lambda_{1}z$ for some nonzero integer $\lambda_{1}$. Comparing
     \[\bar{Q}_{0}z={\textstyle \sum_{i=0}^{n-3}} a_{i}\bar{Q}_{0}S^{i}e+a_{n-2}\bar{Q}_{0}S^{n-2}e={\textstyle \sum_{i=0}^{n-3}}(-1)^{i} a_{i}\ell_{0}S^{i}e+a_{n-2}\ell_{0}S^{n-2}e\] and
     \[\lambda_{1}z={\textstyle \sum_{i=0}^{n-3}} a_{i}\lambda_{1}S^{i}e+a_{n-2}\lambda_{1}S^{n-2}e.\]
     
     Then we get $\lambda_{1}=\ell_{0}$, i.e., $\bar{Q}_{0}z=\ell_{0}z$. Thus, we have
     \begin{center}
         $(\ell_{0}z)^{\rm T}\alpha=(\bar{Q}_{0}z)^{\rm T}\alpha=z^{\rm T}(\bar{Q}_{0}\alpha)=-z^{\rm T}(\ell_{0}\alpha)$,
     \end{center}
     as $Q_{0}$ is a symmetric matrix. Since $p\nmid \ell_{0}$ and $p\nmid 2$, it is easy to see $z^{\rm T}\alpha=0$. It follows that $z$ is orthogonal to  $V=\mathbb{F}_{p}^{n}$; a contradiction. Thus, one obtains $z^{\rm T}z\ne 0$. This completes the proof.
\end{proof}

We are now ready to present the proof of Theorem \ref{l mid l0}.

\begin{proof}[Proof of the Theorem \ref{l mid l0}]
   Suppose to the contrary, there exists an odd prime $p$ that divides $\ell$ but $ p\nmid\ell_{0}$. Let $\bar{Q}=\ell Q$. By $Q^{\rm T}SQ=S(\Delta)$ and $Q^{\rm T}e=e$, we have $W^{\rm T}\bar{Q}=\ell W^{\rm T}(\Delta)\equiv0\pmod{p}$. Let $v$ a vector such that $W^{\rm T}v\equiv0\pmod{p}$. Note that $\mathrm{rank}_{p}(W^{\rm T})=n-1$, then it is easy to see that each column $q_{i}$ of $\bar{Q}$ is linearly dependent with $v$ over $\mathbb{F}_{p} $, and then there exists an $a_{i}\in\mathbb{F}_{p}$ such that $v=a_{i}q_{i}$.
Combining $\bar{Q}^{\rm T}\bar{Q}=\ell^{2}I$, we get
\begin{equation}
    v^{\rm T}v\equiv a_{i}^{2}q_{i}^{\rm T}q_{i}=a_{i}^{2}\ell^{2}\equiv 0\pmod{p}.
\end{equation}
It contracts Lemma \ref{vv=0}. Thus, if $p$ is an odd prime of $\ell$, then $p$ divides $ \ell_{0}$. Since $\Sigma\in \mathcal{G}_{n}$, which shows $d_{n}$ is square-free, it follows that the level $\ell$ of every $Q\in \mathcal{Q}(\Sigma)$ is also square-free by Lemma \ref{l dn}. Clearly, $\ell\mid\ell_{0}$.  This completes the proof.
\end{proof}

\begin{corollary}[\cite{ref4}]
     Let $\Sigma\in \mathcal{G}_{n}$. If $\Sigma$ is self-converse, then $\Sigma$ is $\mathrm{DGSS}$.
\end{corollary}
\begin{proof}
    When $\Sigma$ is self-converse, one obtains that $\ell_{0}=1$. By Theorem~\ref{l mid l0}, it is easy to see that the
level $\ell(Q)$ of every $Q\in \mathcal{Q}(\Sigma)$ is one. It follows that $\Sigma$ is $\mathrm{DGSS}$ by Corollary \ref{l=1}.
\end{proof}

\section{Proofs of Theorems \ref{2t-1} and \ref{Q=QP}}\label{sec4}
In this section, we present the proofs of Theorems \ref{2t-1} and \ref{Q=QP}. First, we need several lemmas below.

\begin{lemma}\label{vector relationship}
     Let $\Sigma \in \mathcal{G}_{n}$ and $Q\in \mathcal{Q}(\Sigma)$ with level $\ell$. Then there exists an integral vector $\alpha\not\equiv 0 \pmod{\ell}$ which is a solution of $W^{\rm T}x\equiv0\pmod{\ell}$ such that for any solution $z$ of the equation, we have $z\equiv k\alpha\pmod{\ell}$ for some integer $k$.
\end{lemma}
\begin{proof}
     By Lemma~\ref{LQ}, let $W^{\rm T}=USV$, where $U$ and $V$ are unimodular matrices, and
      \[S={\rm diag}( \underbrace{1,1,\dots,1}_{\left \lceil \frac{n}{2}\right \rceil},\underbrace{2,2,\dots,2,2b}_{\left \lfloor \frac{n}{2}  \right \rfloor} ),\]
 is the SNF of $W^{\rm T}$. Then the congruence equation $W^{\rm T}x\equiv0\pmod{\ell}$ is equivalent to $SVx\equiv0\pmod{\ell}$. Let $Vx=y$. Note that $\ell\mid d_{n}=2b$ and $\ell$ is odd, one obtains that the general solution of $Sy\equiv0\pmod{\ell}$ is $y=(0,0,\dots,0,m)^{\rm T}$, where $m$ is any integer. Let $\alpha=V^{-1}(0,0,\dots,0,1)^{\rm T}$. It follows that any solution $z$ of the equation is a multiple of $\alpha$ modulo $\ell$, i.e.,
     $z\equiv m\alpha\pmod{\ell}$ for some integer $m$.

\end{proof}

The following lemma lies at the heart of our proof.
\begin{lemma}\label{l=l0/l1}
    Let $\Sigma, \ \Delta,\ \Gamma \in \mathcal{G}_{n}$ be pairwise generalized cospectral. There exist matrices $Q\in\mathcal{Q} (\Sigma)$ with the level $\ell$ and $Q_{1}\in\mathcal{Q} (\Sigma)$ with the level $\ell_{1}$ such that
    $Q^{\rm T}S(\Sigma)Q=S(\Delta),$ and
    $Q_{1}^{\rm T}S(\Sigma)Q_{1}=S( \Gamma )$ , where $\ell_{1} \mid \ell $.
    Let $Q_{2}=Q_{1}^{\rm T}Q$. Then we
     have $Q_{2}^{\rm T}S(\Gamma)Q_{2}=S(\Delta)$. Moreover, the level of $Q_{2}$ is $\ell/\ell_{1}$.
\end{lemma}

\begin{proof}

Note that $Q_{2}=Q_{1}^{\rm T}Q$. It follows that
    \begin{align*}
       Q_{2}^{\rm T}S(\Gamma) Q_{2}&=(Q_{1}^{\rm T}Q)^{\rm T}S(\Gamma)(Q_{1}^{\rm T}Q)\\
        &=Q^{\rm T}Q_{1}S(\Gamma)Q_{1}^{\rm T}Q\\
        &=Q^{\rm T}S(\Sigma)Q=S(\Delta).
        \end{align*}

Suppose that the level of $Q_{2}$ is $\tilde{\ell}$. On the one hand, left-multiplying both sides of the equation
$Q_{2}=Q_{1}^{\rm T}Q$ by the matrix $Q_{1}$, we obtain $Q_{1}Q_{2}=Q$ as $Q_{1}$ is orthogonal. It is easy to see that $\ell\mid \ell_{1}\tilde{\ell} $ and hence
 $(\ell/\ell_{1})\mid \tilde{\ell}$.

On the other hand,
\[\frac{\ell}{\ell_{1}} Q_{2}=\frac{\ell}{\ell_{1}}Q_{1}^{\rm T}Q=\frac{(\ell_{1} Q_{1})^{\rm T}(\ell Q)}{\ell_{1}^{2}} .\]
By the above equation, to prove  $\tilde{\ell}\mid \ell/\ell_{1}$, it suffices to prove the matrix $(\ell_{1} Q_{1})^{\rm T}(\ell Q)/{\ell_{1}^{2}}$ is an integral matrix. Let $\ell_{1}Q_{1}=\left [\beta_{1},  \beta_{2},\dots,\beta_{n}\right ] $ and $\ell Q=\left [\gamma  _{1},  \gamma  _{2},\dots,\gamma  _{n}\right ]$.\\

\noindent
    \textbf{ Claim:} $\beta_{i}^{\rm T}\gamma _{j}\equiv 0 \pmod{\ell^{2}_{1}}$, for $1\le i,j\le n$.\\

Since $\ell_{1} \mid \ell$, we get $W^{\rm T}(\ell Q)=\ell W^{\rm T}(S(\Delta))\equiv0\pmod{\ell_{1}}$. Let $\alpha $ be a solution of $W^{\rm T}x\equiv0\pmod {\ell_{1}}$. By Lemma \ref{vector relationship}, we have $\beta_{i}\equiv k_{i}\alpha\pmod{\ell_{1}}$ for some integer $k_{i}$, and $\gamma_{j}\equiv m_{j}\alpha\pmod{\ell_{1}}$ for some integer $m_{j}$. Then
\[\ell_{1}\mid \beta_{i}-k_{i}\alpha,\ \ell_{1}\mid \gamma_{j}-m_{j}\alpha.\]
 In what follows, let $p$ be any prime factor of $\ell_{1}$. Then we proceed by considering the following four cases regarding whether $k_{i}$ and $m_{j}$ are zero or not over $\mathbb{F}_{p}$.\\

 \noindent
\textbf{Case 1.} $p\nmid k_{i}$ and $p\nmid m_{j}.$
Under the condition, we have
$ p\mid m_{j}\beta_{i}-k_{i}\gamma_{j}$ by eliminating $\alpha$,
and hence
\begin{equation}\label{br}
    p^{2}\mid (m_{j}\beta_{i}-k_{i}\gamma_{j})^{\rm T}(m_{j}\beta_{i}-k_{i}\gamma_{j})=m_{j}^{2}\beta_{i}^{\rm T}\beta_{i}+2m_{j}k_{i}\beta_{i}^{\rm T}\gamma_{j}+k_{i}^{2}\gamma_{j}^{\rm T}\gamma_{j}.
\end{equation}
Since $Q$ and $Q_{1}$ are orthogonal, then $p^{2}\mid  \beta_{i}^{\rm T} \beta_{i} $ and $p^{2}\mid\gamma_{i}^{\rm T} \gamma_{i} $, it follows from Eq.~(\ref{br}) that
$ p^{2}\mid \beta_{i}^{\rm T}\gamma_{j}. $\\

\noindent
\textbf{Case 2.} $p\mid k_{i}$ and $p\nmid m_{j}.$ Let $\beta_{s}$ be the $s$-th column of $\ell_{1} Q_{1}$ such
that $\beta_{s}\not\equiv0\pmod{p}$. Then $ \beta_{s}\equiv k_{s}\alpha\pmod{p}$
 for some integer $p\nmid k_{s}$. Indeed, the existence of such a vector is necessary; otherwise, $\ell_{1} Q_{1}/p$ is an integral matrix, which is a contradiction. \\
 By eliminating $\alpha$, we derive  $p\mid \beta_{i}-k_{s}\gamma_{j}+m_{j}\beta_{s}.$ Then
\begin{equation}
     p^{2}\mid (\beta_{i}-k_{s}\gamma_{j}+m_{j}\beta_{s})^{\rm T}(\beta_{i}-k_{s}\gamma_{j}+m_{j}\beta_{s}).
\end{equation}
and hence,
\begin{equation}
     p^{2}\mid \beta_{i}^{\rm T}\beta_{i}+k_{s}^{2}\gamma_{j}^{\rm T}\gamma_{j}+m_{j}^{2}\beta_{s}^{\rm T}\beta_{s}-2k_{s}\beta_{i}^{\rm T}\gamma_{j}+2m_{j}\beta_{i}^{\rm T}\beta_{s}-2m_{j}k_{s}\beta_{s}^{\rm T}\gamma_{j}.
\end{equation}
As established in Case 1, we have $p^{2}\mid\beta_{s}^{\rm T}\gamma_{j}$. Since $Q$ is orthogonal, it shows that $\beta_{i}^{\rm T}\beta_{s}=0$. Note that $p^{2}\mid \beta_{s}^{\rm T} \beta_{s} $.
 Therefore, we obtain
 $p^{2}\mid \beta_{i}^{\rm T}\gamma_{j}.$\\

 \noindent
\textbf{Case 3.} $p\nmid k_{i}$ and $p\mid m_{j}.$ The proof of this case is similar to Case 2. The result $p^{2}\mid \beta_{i}^{\rm T}\gamma_{j}$ also hold in this case.\\

\noindent
\textbf{Case 4.} $p\mid k_{i}$ and $p\mid m_{j} .$ In this case, $p\mid \beta_{i}$ and $p\mid \gamma _{j}$,
 it is easy to see $p^{2}\mid \beta_{i}^{\rm T}\gamma _{j}.$\\

Therefore, for any odd prime factor $p$ of $\ell_1$, one obtains that $\beta_{i}^{\rm T}\gamma _{j}\equiv 0 \pmod{p^{2}},$  for $1\le i,j\le n$. It follows immediately that $\beta_{i}^{\rm T}\gamma _{j}\equiv 0 \pmod{\ell_{1}^{2}}$. The claim follows.

Futhermore, it implies that the $(\ell/\ell_{1})Q_{2}=(\ell_{1} Q_{1})^{\rm T}(\ell Q)/{\ell_{1}^{2}}$ is integral matrix, and hence $\tilde{\ell}\mid \ell/\ell_{1} $. It follows $\tilde{\ell}=\ell/\ell_{1}$. The proof is complete.
\end{proof}

\begin{corollary}\label{Q0=Q1Q2}
     Let $\Sigma, \ \Delta\in \mathcal{G}_{n}$ and $Q_{0}$ be defined as in \eqref{EE1} with $\ell_{0}$. Suppose that there exists a matrix $Q_{1}\in \mathcal{Q}(\Sigma)$ with level $\ell_{1}$ such that
    $Q_{1}^{\rm T}S(\Sigma)Q_{1}=S(\Delta).$ Let $Q_{2}=Q_{1}^{\rm T}Q_{0}$.  Then we have $Q_{2}S(\Sigma)Q_{2}^{\rm T}=S(\Delta^{\rm T})$. Moreover, the level of $Q_{2}$ is $\ell_{0}/\ell_{1}$.
\end{corollary}
\begin{proof}
    By Lemma \ref{l=l0/l1}, since $Q_{0}^{\rm T}S(\Sigma)Q_{0}=S(\Sigma^{\rm T})$, we have $Q_{2}^{\rm T}S(\Delta)Q_{2}=S(\Sigma^{\rm T})$. Thus, $Q_{2}S(\Sigma)Q_{2}^{\rm T}=S(\Delta^{\rm T})$ as $S^{\rm T}=-S$. The second assertion can be obtained from the above Lemma immediately.
\end{proof}

Now, we are ready to present the proof of Theorem~\ref{Q=QP}.
    \begin{proof}[Proof of Theorem~\ref{Q=QP}]
       Since $\Delta,\ \Gamma \in \mathcal{G}_{n}$ are pairwise generalized cospectral, there exists a regular orthogonal rational matrix $Q_{2}$ such that
       \begin{equation}\label{11}
            Q_{2}^{\rm T}S(\Gamma)Q_{2}=S(\Delta).
       \end{equation}
       According to Lemma~\ref{QSQ}, the regular orthogonal rational matrix $Q_1$ satisfying the Eq.~(\ref{11}) is unique. Thus, we obtain $Q_{2}=Q_{1}^{\rm T}Q$ and the level of $Q_{2}$ is $\ell/\ell=1$ by Lemma \ref{l=l0/l1}. Hence, $Q_{2}$ is a permutation matrix, which implies that $ \Delta$ is isomorphic to $\Gamma $.
\end{proof}

Combining Theorem \ref{l mid l0} and Theorem~\ref{Q=QP}, we finally present the proof of Theorem \ref{2t-1}.
\begin{proof}[Proof of Theorem \ref{2t-1}]
    Let $\Delta$ be any oriented graph that is generalized cospectral
with $\Sigma$. Then there exists a $Q\in \mathcal{Q}(\Sigma )$ with level $\ell=\ell(Q)$ such that $Q^{\rm T}S(\Sigma)Q = S(\Delta)$.
     According to Theorem \ref{l mid l0}, the level $\ell$ of $Q$ divides $\ell_{0}=p_{1}p_{2}\cdots p_{t}$. Hence, there are at most $2^{t}-1$ possible choices for the level $\ell$. While Theorem~\ref{Q=QP} indicates that, for any two matrices $Q_{1},\ Q_{2}\in\mathcal{Q}(\Sigma )$, if $\ell(Q_{1})=\ell(Q_{2 })=\ell$, then the two oriented graphs with
adjacency matrices $Q_{1}^{\rm T}S(\Sigma)Q_{1}$ and $Q_{2}^{\rm T}S(\Sigma)Q_{2}$, respectively, are actually isomorphic.  Therefore, up to isomorphism, the number of generalized cospectral mates of $\Sigma$ is at most $2^{t}-1$. This completes the proof.
   \end{proof}
    As an immediate consequence, we have
   \begin{corollary}
       If $t=1$, or equivalently $\ell_0$ is a prime, then $\Sigma\in\mathcal{G}_{n}$ is $\mathrm{WDGSS}$.
   \end{corollary}
\begin{proof}
    By Theorem~\ref{2t-1}, when $t=1$, the number of cospectral mates
of $\Sigma$ is at most 1. If $\Sigma$ is self-converse, then it is $\mathrm{DGSS}$ and
hence $\mathrm{WDGSS}$. If $\Sigma$ is non-self-converse, it has a generalized
cospectral mate $\Sigma^{\rm T}$. Thus, if the $\ell_0$ is an odd prime,
then $\Sigma$ is $\mathrm{WDGSS}$.
\end{proof}

\section{A criterion for $\mathrm{WDGSS}$ oriented graphs}\label{sec5}
Based on the previous analysis, in this section, we are able to give a partial answer to Problem~\ref{P1}. Specifically, we give an efficient criterion for determining whether an oriented graph $\Sigma \in \mathcal{G}_{n}$ with level $\ell_{0}=p_1p_2\cdots p_t$ is $\mathrm{WDGSS}$, and finding out all its generalized cospectral mates if it is not, for $t\leq 3$.

Recall that in Corollary~\ref{Q0=Q1Q2}, we know that whenever there is a cospectral mate $\Delta$ of $\Sigma$ and
$Q_{1}^{\rm T}S(\Sigma)Q_{1}=S(\Delta)$ for a regular rational orthogonal matrix $Q_1$ with
level $\ell_1$. Then $Q_0$ can be written as the product of two regular
rational orthogonal matrices $Q_1$ and $Q_2=Q_1^{\rm T}Q_0$ with levels $\ell$ and $\ell_0/\ell$, respectively. In other words, if the matrix $Q_0$ cannot be decomposed into the product of two matrices with smaller levels, then $\Sigma$ has no
generalized cospectral mates other than its converse $\Sigma^{\rm T}$.
\begin{definition}
   A regular rational orthogonal matrix $Q$ of level $\ell$ is called \emph{factorable}
    if it can be written as the product of two regular rational orthogonal matrices of smaller levels, and neither of which is a permutation matrix.
\end{definition}

\begin{theorem}
     Let $\Sigma\in\mathcal{G}_{n}$ and $Q_0$ be defined as in Eq.~(\ref{EE1}). If $Q_{0}$ is not factorable, then the oriented graph $\Sigma$ is $\mathrm{WDGSS}$
\end{theorem}
\begin{proof}

    Suppose to contrary that $\Sigma\in\mathcal{G}_{n}$ is not $\mathrm{WDGSS}$. Let $\Delta$ be another generalized cospectral
mates of $\Sigma$ other than $\Sigma^{\rm T}$. Then there exists a $Q_{1}\in \mathcal{Q} (\Sigma)$ with the level $\ell_{1}\ne\ell_{0}$ and $\ell_{1}\ne1$ such that $Q_{1}^{\rm T}S(\Sigma)Q_{1}=S(\Delta)$. By Corollary \ref{Q0=Q1Q2}, we know that there exists another graph $\Delta^{\rm T}$ with $Q_{2}SQ_{2}^{\rm T}=S(\Delta^{\rm T})$, where $Q_{2}=Q_{1}^{\rm T}Q_{0}\in \mathcal{Q} (\Sigma)$ with the level $\ell_{2}=\ell_{0}/\ell_{1}\ne1$. Moreover, we have $Q_0=Q_1Q_2$. It contradicts the assumption that $Q_{0}$ is not factorable.
    Thus, $\Sigma\in\mathcal{G}_{n}$ is $\mathrm{WDGSS}$.
\end{proof}
\begin{remark}
    An efficient algorithm was given in \cite{ref3} to determine
whether $Q_0$ is factorable or not, for $t \le 3$, and it can find out all
generalized cospectral mates of $\Sigma$ when it is not $\mathrm{WDGSS}$.

\end{remark}

    In what follows, we consider a special case that $\ell_0=p_1p_2$. We shall describe the relationship between all the cospectral mates of $\Sigma$, when $\Sigma$ is not $\mathrm{WDGSS}$.
    
     Suppose $\Sigma$ is not  $\mathrm{WDGSS}$, and $\Sigma^{\rm T}$, $\Delta$ and $\Delta^{\rm T}$ are all its generalized cospectral mates. Then by Lemma~\ref{QSQ}, there exists a matrix $Q_{1}\in\mathcal{Q}(\Sigma )$ with level $p_1$ such that $Q_{1}^{\rm T}S(\Sigma)Q_{1}=S(\Delta)$, and there exists a matrix $Q_{2}\in\mathcal{Q}(\Sigma )$ with level $p_2$ such that $ Q_{2}^{\rm T}S(\Delta)Q_{2}=S(\Sigma^{\rm T})$.
    
    \begin{lemma}\label{Q0_new}
     Suppose that $\hat{Q}_{0}$ is a regular rational orthogonal matrix such that $\hat{Q}_{0}^{\rm T}S(\Delta^{\rm T})\hat{Q}_{0}=S(\Delta)$. Then $\hat{Q}_{0}=Q_2Q_{1}$ and $\ell(\hat{Q}_{0})=p_{1}p_{2}$.
    \end{lemma}
\begin{proof}
    Let $Q=Q_2Q_{1}$.
 It is obvious that $Q$ is a rational orthogonal and regular matrix. Then
     \begin{align*}
        Q^{\rm T}S(\Delta^{\rm T})Q&=(Q_2Q_{1})^{\rm T}S(\Delta^{\rm T})(Q_2Q_{1})\\
        &=Q_1^{\rm T}Q_2^{\rm T}S(\Delta^{\rm T})Q_2Q_{1}\\
        &=Q_{1}^{\rm T}S(\Sigma)Q_{1}=S(\Delta)
    \end{align*}
    By Lemma \ref{QSQ}, where $Q$ such that $Q^{\rm T}S(\Delta^{\rm T})Q=S(\Delta)$ is unique, then $\hat{Q}_{0}=Q_2Q_{1}$.

    Let the level of $\hat{Q}_{0}$ be $\hat{\ell}$. Note that $\ell(Q_{1})=p_{1}$ and $\ell(Q_{2})=p_{2}$, it follows $\hat{\ell}\mid p_{1}p_{2}$.\\
   \textbf{Claim:} $\hat{\ell}= p_{1}p_{2}$.

   Suppose that $\hat{\ell}= p_{2}$, then $p_{2}\hat{Q}_{0}=p_{2}Q_{2}Q_{1}$ is an integral matrix, which implies that $p_{1}p_{2}\hat{Q}_{0}=p_{1}p_{2}Q_{2}Q_{1}$ is zero matrix over $\mathbb{F}_{p_{1}}$. Then we get
   \begin{equation}
      p_{1}p_{2}Q_{2}Q_{1}\equiv (p_{2}Q_{2})(p_{1}Q_{1})\equiv0\pmod{p_{1}}.
   \end{equation}
   As $p_{1}$ and $p_{2}$ are relatively prime, we have $\mathrm{rank}_{p_{1}}(p_{2}Q_{2})=n$. Thus, the solution space of $(p_{2}Q_{2})x\equiv0\pmod{p_{1}}$ has only zero solution. It follows that any column $q_{i}$ of $p_{1}Q_{1}$ is $q_{i}\equiv0\pmod{p_{1}}$, and hence we get $\mathrm{rank}_{p_{1}}(p_{1}Q_{1}))=0$, which contradicts that the level of $Q_{1}$ is $p_{1}$. It shows the assumption $\hat{\ell}= p_{2}$ does not hold. Likewise, we get $\hat{\ell}\ne p_{1}$ in a similar way. By Theorem \ref{DGSS}, we know that when $\Delta\in \mathcal{G}_{n}$ is self-converse, then $\Delta$ does not have cospectral mate. Therefore, the level of $\hat{Q}_{0}$ cannot be equal to 1, and hence $\hat{\ell}=p_{1}p_{2}$.
 \end{proof}

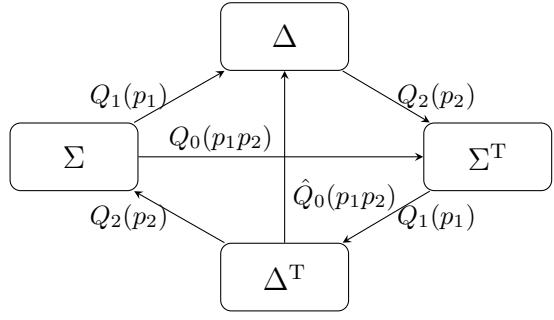
\begin{figure}[ht]\label{fig}
    \centering
  \thispagestyle{empty}

\tikzstyle{startstop} = [rectangle, rounded corners, minimum width = 1.7cm, minimum height=0.9cm,text centered, draw = black]
\tikzstyle{io} = [rectangle, rounded corners, minimum width = 1.7cm, minimum height=0.9cm,text centered, draw = black]
\tikzstyle{process} = [rectangle, rounded corners, minimum width = 1.7cm, minimum height=0.9cm,text centered, draw = black]
\tikzstyle{decision} = [rectangle, rounded corners, minimum width = 1.7cm, minimum height=0.9cm,text centered, draw = black]
\tikzstyle{arrow} = [->,>=stealth]
\begin{tikzpicture}[node distance=2cm]
\node[startstop](start){$\Sigma$};
\node[io, right of = start, xshift = 3.5cm](in1){$\Sigma^{\rm T}$};
\node[process, xshift = 2.8cm ,yshift = 1.6cm](pro1){$\Delta$};
\node[decision, xshift = 2.8cm, yshift = -1.6cm](dec1){$\Delta^{\rm T}$};

\draw [arrow] (start)-- node [xshift = -0.8cm ,yshift = 0.25cm] {\footnotesize{$Q_0(p_1p_2)$}} (in1);
\draw [arrow] (dec1) -- node [left] {\footnotesize{$Q_2(p_2)$}}(start);
\draw [arrow] (start) --
node [left] {\footnotesize{$Q_1(p_1)$}}(pro1);
\draw [arrow] (pro1) --
node [right] {\footnotesize{$Q_2(p_2)$}}(in1);
\draw [arrow] (in1) -- node [right] {\footnotesize{$Q_1(p_1)$}} (dec1);
\draw [arrow] (dec1) -- node [xshift = 0.8cm ,yshift = -0.5cm] {\footnotesize{$\hat{Q}_0 (p_1p_2)$}} (pro1);

   \end{tikzpicture}
   \caption{\label{fig:Relationship} The relationship between all cospectral mates of $\Sigma$ with $\ell_0=p_1p_2$.}
\end{figure}

     Let $\Sigma$, $\Sigma^{\rm T}$, $\Delta$ and $\Delta^{\rm T}$ be pairwise generalized cospectral. Then there exists a matrix $Q_{0}\in\mathcal{Q}(\Sigma )$ with level $p_1p_2$ such that
     \[Q_0^{\rm T}S(\Sigma)Q_0=S(\Sigma^{\rm T}).\]
     Since $\Sigma$ and $\Delta$ be  generalized cospectral,
     then there exists a matrix $Q_{1}\in\mathcal{Q}(\Sigma )$ with level $p_1$ such that
    \[Q_{1}^{\rm T}S(\Sigma)Q_{1}=S(\Delta).\]
    Transposing both sides of the equation, we can obtain $Q_{1}^{\rm T}S(\Sigma^{\rm T})Q_{1}=S(\Delta^{\rm T}).$
    By Corollary~\ref{Q0=Q1Q2}, we know that there exists a matrix $Q_2=Q_1^{\rm T}Q_0$ with level $p_2$ such that
    \[Q_{2}^{\rm T}S(\Delta^{\rm T})Q_{2}=S(\Sigma).\]
    Likewise, transposing both sides of the equation, we have $Q_{2}^{\rm T}S(\Delta)Q_{2}=S(\Sigma^{\rm T}).$
    It follows from Lemma~\ref{Q0_new} that there exists a matrix $\hat{Q}_0=Q_2Q_1$ with the level $p_1p_2$ such that
    \[\hat{Q}_0^{\rm T}S(\Delta^{\rm T})\hat{Q}_0=S(\Delta).\]
     In summary, the relationship between $\Sigma$, $\Sigma^{\rm T}$, $\Delta$ and $\Delta^{\rm T}$ is shown in Fig.~\ref{fig:Relationship}.

     Next, we shall describe a method for identifying  all cospectral mates of $\Sigma$ and further elucidate the relationships among these cospectral mates. We shall focus on the special case that $\ell_{0}=p_{1}p_{2}$, where $p_{1}$ and $p_{2}$ are distinct odd primes.
     The method also applies to the case that $\ell_0$ factors into three distinct odd primes, i.e., $\ell_{0}=p_{1}p_{2}p_3$, but fails for $\ell_{0}=p_{1}p_{2}\cdots p_t$ for $t>3$.
     
    We need the following definitions.
\begin{definition}[\cite{ref3}]
   {\rm A regular rational orthogonal matrix $Q$ of level $p$ is called a \emph{primitive matrix}
if $\mathrm{rank}_{p}(pQ)=1$, where $p$ is a prime.}
\end{definition}
\begin{definition}[\cite{ref3}]
    Let $v\not\equiv0~({\rm mod}~p)$ be an $n$-dimensional integral vector and $Q$ be a primitive matrix. We say $Q$ can be \emph{generated} from $v$ (or $v$ \emph{generates} $Q$) if each column of $pQ$ is a multiple of $v$ over $\mathbb{F}_{p}$.
\end{definition}

    Now, suppose there exists a regular rational orthogonal matrix $Q_1$ with level $\ell_1=p_1$ such that $Q_1^{\rm T}S(\Sigma)Q_1=S(\Delta)$.
    Then it follows that $Q_1^{\rm T}W(\Sigma)=W(\Delta)$, or equivalently, \[W(\Sigma)^{\rm T}(p_1Q_1)=p_1W(\Delta)^{\rm T}\equiv 0~({\rm mod}~p_1).\]

Note that ${\rm rank}_{p_1}W(\Sigma)=n-1$, the nullspace of $W(\Sigma)^{\rm T}$ is 1-dimensional, over $\mathbb{F}_{p_1}$. It follows that $Q_1$ is a primitive rational orthogonal matrix, and the columns of $p_1Q_1$ lie in this space.

Conversely, suppose that the nullspace of $W(\Sigma)^{\rm T}$ is spanned by a vector $v$. Then the algorithm in \cite{ref3} can be applied to
determine whether the vector $v$ can generate a primitive rational orthogonal matrix $Q_1$. If the answer is ``NO", then $Q_0$ is not factorable and $\Sigma$ is $\mathrm{WDGSS}$.
If the answer is ``YES", then we have to check whether $Q_1^{\rm T}S(\Sigma)Q_1$ is the skew-adjacency matrix of an oriented graph, say $\Delta$. If the answer is ``YES", let $Q_2=Q_1^{\rm T}Q_0$, then $Q_2S(\Sigma)Q_2^{\rm T}=S(\Delta^{\rm T})$. We have found all generalized cospectral mates of $\Sigma$, namely, $\Delta$, $\Delta^{\rm T}$ and $\Sigma^{\rm T}$. Otherwise, $\Sigma$ is $\mathrm{WDGSS}$.

\section{Some examples}\label{sec6}
In this section, we give some examples for illustrations. In particular, these examples show that the upper bound in Theorem~\ref{2t-1} is attainable.
We would like to mention that the computation of all the regular rational orthogonal matrices involved (except for $Q_0$) is based on the algorithm proposed in \cite{ref3}.

\begin{exam}
   Let $\Sigma$ be an oriented graph of order $n=6$ given as in Fig~\ref{fig:l0=1}.
    \begin{figure}[ht]
        \centering
\includegraphics[width=0.5\linewidth]{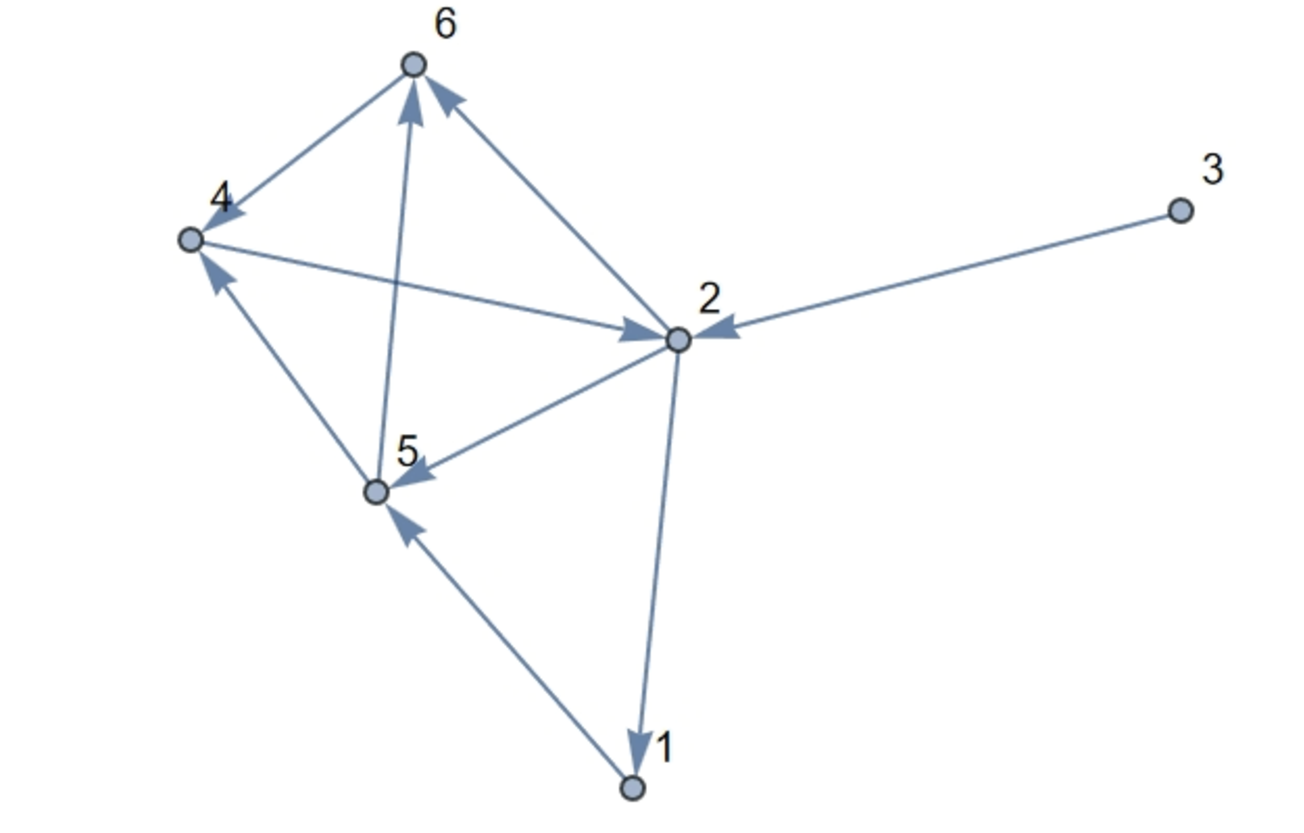}
        \caption{A non-self-converse oriented graph $\Sigma$ which is $\mathrm{WDGSS}$.}
        \label{fig:l0=1}
    \end{figure}\\
    The skew-adjacency matrix $S(\Sigma)$ of
$\Sigma$ is given as follows:
\[S(\Sigma)=\begin{bmatrix} 0& -1& 0& 0& 1& 0\\
       1& 0& -1& -1& 1& 1\\
       0& 1& 0& 0& 0& 0\\
       0& 1& 0& 0& -1& -1\\
   -1& -1& 0& 1& 0&1\\
   0&-1& 0& 1& -1& 0\\
   \end{bmatrix}.\]
It can be computed that the $\mathrm{SNF}$ of $W(\Sigma)$ is
       \[
           \mathrm{diag} \left ( 1,1,1,2,2,2\times 523\right ),
       \]
       and the $Q_{0}$ such that $Q_{0}^{\rm T}S(\Sigma)Q_0=S(\Sigma^{\rm T})$ is given as follows:
          \[ Q_{0}=\frac{1}{523} \begin{bmatrix}
359& -198& 146& 48& 268& -100\\ -198& -137& 138& 160&196& 364\\ 146& 138& -181& 442& -60& 38\\
48& 160& 442&  88& -206& -9\\
268& 196& -60& -206& -17& 342\\ -100& 364& 38& -9&
  342& -112\\
           \end{bmatrix}.\]
   It is easy to see that the level of $Q_0$ is $\ell_{0} = 523$, which is a prime.
 Thus, $\Sigma$ is $\mathrm{WDGSS} $.

\end{exam}

\begin{exam}
Let $\Sigma$ be an oriented graph of order $n=6$ given as in Fig~\ref{fig:l0=2}, and the skew-adjacency matrix $S(\Sigma)$ of
$\Sigma$ be given as follows:
\begin{figure}[ht]
    \centering
\includegraphics[width=0.5\linewidth]{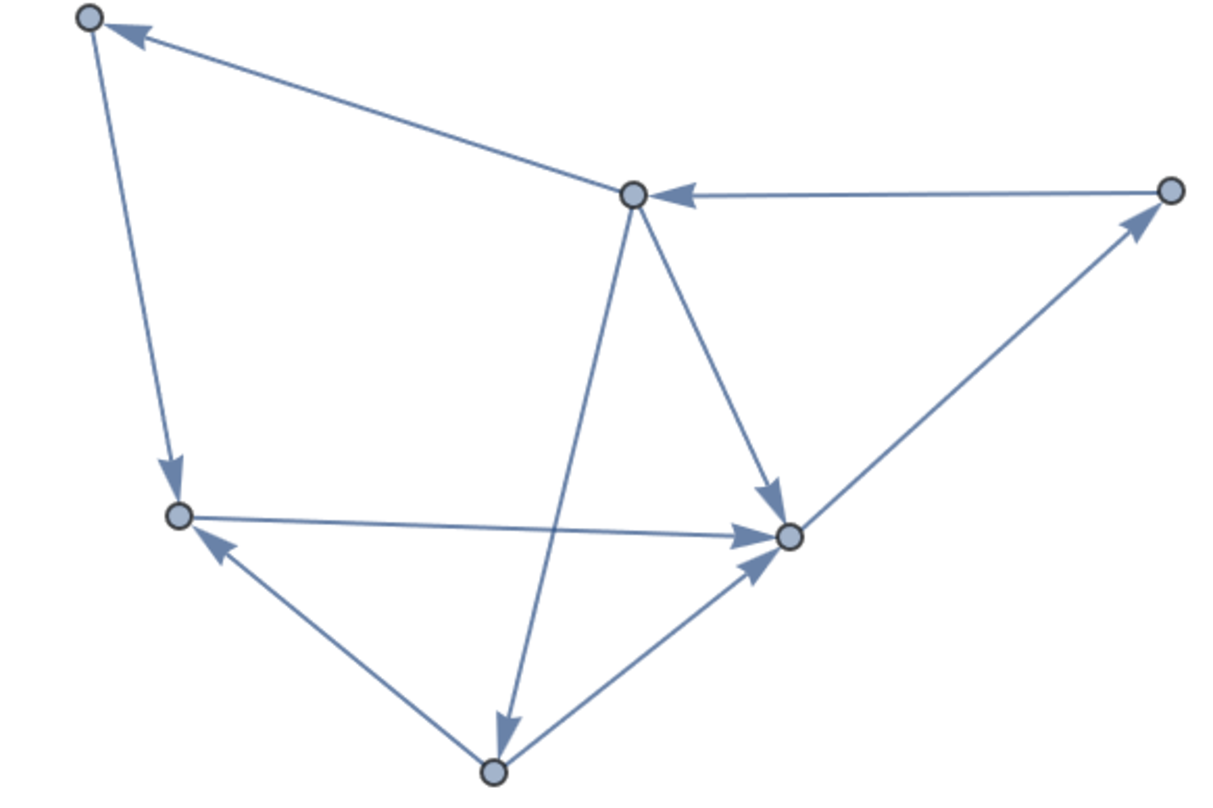}
    \caption{A non-self-converse oriented graphs $\Sigma$.
}
    \label{fig:l0=2}
\end{figure}

\[
S(\Sigma)=\begin{bmatrix}
       0& 1& -1& 0& 0& 0\\
-1& 0& 0& 1& -1& 0\\
1& 0&0&1&  1& -1\\
0& -1& -1& 0& -1& 1\\
0&1&-1& 1& 0&0\\
0&0&1& -1& 0&0\\
\end{bmatrix}.
\]
It can be computed that the $\mathrm{SNF}$ of $W(\Sigma)$ is
       \[
           \mathrm{diag} \left ( 1,1,1,2,2,2\times 3\times7\right ),
       \]
       and the $Q_{0}$ is
      \[
           Q_{0}=\frac{1}{21} \begin{bmatrix}
2& -9& -1& 11& 15& 3\\
           -9&  9&  15&  3&  6&  -3\\
           -1&  15&  -10&  5&  3& 9\\
           11&  3&  5&  -13&  9&  6\\
           15&  6&  3&  9&  -3&  -9\\
           3&  -3&  9& 6&  -9&  15\\
           \end{bmatrix}.\]
    It is easy to see that $\ell_{0} = 3\times7$.
    Let 
    \[
    Q_{1}=\frac{1}{3}\begin{bmatrix}
2& 2& -1&0& 0&0\\
0& 0& 0& 3& 0& 0\\
2& -1& 2& 0& 0& 0\\
-1& 2& 2& 0& 0& 0\\
0& 0& 0& 0& 3& 0\\
0& 0& 0& 0& 0& 3\\
\end{bmatrix}~{\rm and}~
Q_2=\frac{1}{7}\begin{bmatrix}
 -1 & 1 & -3 & 5 & 3 & 2 \\
 3 & -3 & 2 & -1 & 5 & 1 \\
 2 & 5 & -1 & -3 & 1 & 3 \\
 -3 & 3 & 5 & 1 & 2 & -1 \\
 5 & 2 & 1 & 3 & -1 & -3 \\
 1 & -1 & 3 & 2 & -3 & 5 \\
\end{bmatrix}.\]
   Then it is easy to verify that $Q_0=Q_1Q_2=Q_2^{\rm T}Q_1^{\rm T}$. Moreover, we have
 \[
S(\Delta_{1})=Q_{1}^{\rm T}S(\Sigma)Q_{1}=\begin{bmatrix}
0& 1& 0& 1& 1& -1\\
-1& 0& -1& 0& -1& 1\\
0& 1& 0& -1& 0& 0\\
-1& 0& 1& 0& -1& 0\\
-1& 1& 0& 1& 0& 0\\
1& -1& 0& 0& 0& 0\\
\end{bmatrix} .
\]
\[
S(\Delta_{1}^{\rm T})=Q_{2}S(\Sigma)Q_{2}^{\rm T}=\begin{bmatrix}
0& -1& 0& -1& -1& 1\\
1& 0& 1& 0& 1& -1\\
0& -1& 0& 1& 0& 0\\
1& 0& -1& 0& 1& 0\\
1& -1& 0& -1& 0& 0\\
-1& 1& 0& 0& 0& 0\\
\end{bmatrix}.\]
Thus, all the cospectral mates of $\Sigma$ are $\Delta_1$, $\Delta_1^{\rm T}$ and $\Sigma^{\rm T}$. There are $2^2-1=3$ in total.
This example illustrates the relationship shown in the Fig~\ref{fig:Relationship}.
\end{exam}

\begin{exam}
       Let $\Sigma$ be an oriented graph of order $n = 7$ given as in Fig~\ref{fig:l0=3}.
      \begin{figure}[ht]
          \centering
          \includegraphics[width=0.5\linewidth]{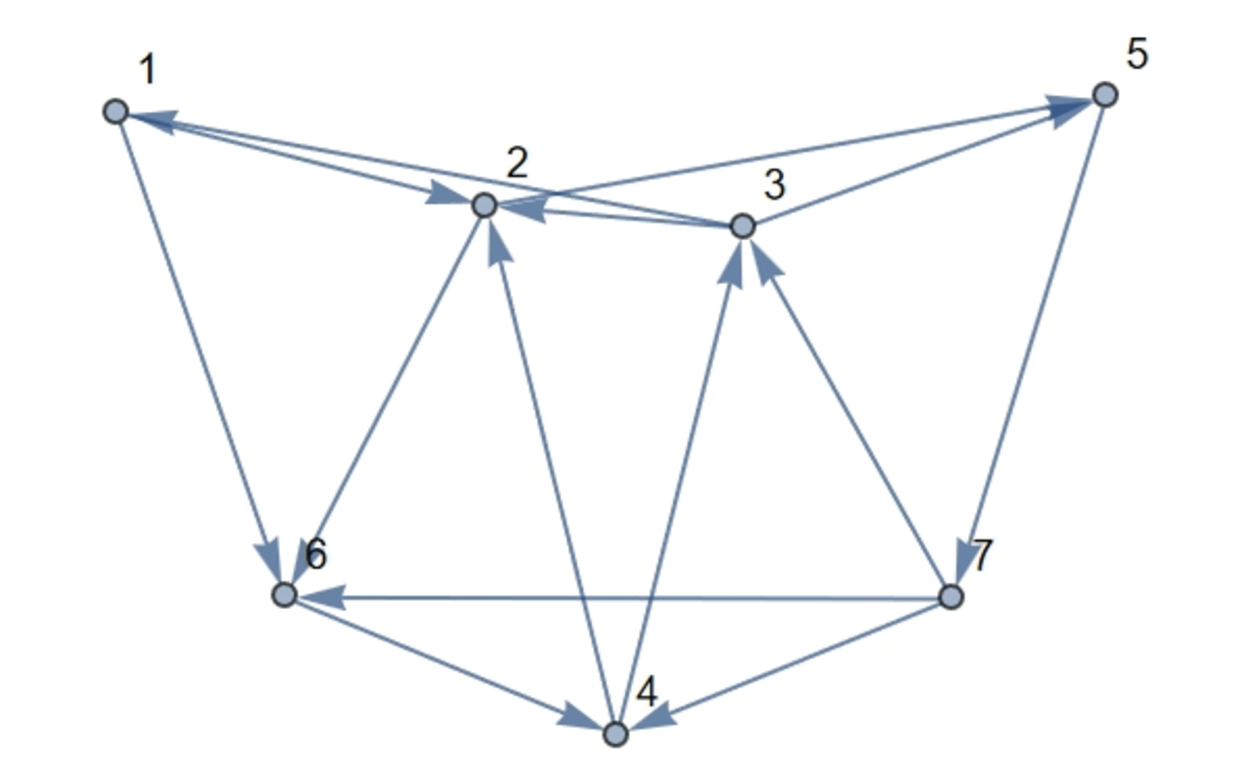}
          \caption{A non-self-converse oriented graphs $\Sigma$}
          \label{fig:l0=3}
      \end{figure}\\
       The skew-adjacency matrix $S(\Sigma)$ of $\Sigma$ is given as follows:
       \begin{center}
           $S(\Sigma)=\begin{bmatrix}
                    0& 1& -1& 0&0& 1& 0\\
       -1& 0& -1& -1&1& 1& 0\\
       1& 1&0& -1& 1& 0& -1\\
       0& 1&1& 0& 0& -1&-1\\
       0& -1& -1& 0& 0& 0& 1\\
       -1& -1& 0& 1& 0& 0&-1\\
       0& 0& 1& 1& -1& 1&0\\
           \end{bmatrix}$.
       \end{center}
 It can be computed that the $\mathrm{SNF}$ of $W$ is
       \begin{center}
          $ \mathrm{diag} \left ( 1,1,1,1,2,2,2\times 7\times11\times13\times211\right )$.
       \end{center}
    The regular rational orthogonal matrix
     \[           Q_{0}=\frac{1}{19201} \begin{bmatrix}
         -4306 & 5754&13736& -3474&-1463&10562& -1608\\
           5754&12502& 3129&9030& -6776& -6230& 1792\\
           13736& 3129& 542& -2670& 9786& 2480&-7802\\
-3474& 9030& -2670& -4631& 10458& -1164& 11652\\
-1463&-6776& 9786& 10458& 9478& -4536& 2254\\
10562& -6230&2480& -1164&-4536& 5396&12693\\
-1608& 1792& -7802& 11652& 2254& 12693& 220\\

           \end{bmatrix}.\]

     It is easy to see the level $\ell_{0}=7\times13\times211$. For $p=11$,  Solving the congruence equation $W^{\rm T}z\equiv0\pmod{11}$ gives $z=(0,6,6,7,2,0,1)^{\rm T}.$
     Note that $z^{\rm T}z\not\equiv0\pmod{11}$. Thus $11\nmid \ell(Q)$ for any $Q\in{\mathcal{Q}(\Sigma)}$. This is precisely predicted by Theorem~\ref{l mid l0} and Lemma~\ref{vv=0}.
     
     For $\ell=7$, let
    \[
    Q_{1}=\frac{1}{7}\begin{bmatrix}
 6 & -1 & 2 & 2 & -2 & 0 & 0 \\
 0 & 0 & 0 & 0 & 0 & 7 & 0 \\
 2 & 2 & -4 & 3 & 4 & 0 & 0 \\
 -2 & -2 & 4 & 4 & 3 & 0 & 0 \\
 0 & 0 & 0 & 0 & 0 & 0 & 7 \\
 -1 & 6 & 2 & 2 & -2 & 0 & 0 \\
 2 & 2 & 3 & -4 & 4 & 0 & 0 \\
\end{bmatrix}.\]

     It is easy to compute that                   
       \begin{center}
$S(\Delta_{1})=Q_1^{\rm T}S(\Sigma)Q_1=\begin{bmatrix}
0& 1& 1& 0&-1&1& 0\\
-1& 0& 0& 1& 0&-1& 0\\
-1& 0& 0& 1& 1& 0& -1\\
0& -1& -1& 0& -1&1& 1\\
1& 0& -1& 1&0&1& 0\\
-1& 1& 0& -1&-1& 0& 1\\
0& 0& 1& -1&0& -1& 0\\
\end{bmatrix}$.
       \end{center}
    By Corollary~\ref{Q0=Q1Q2}, let $\tilde{Q}_{1}:=Q_{0}^{\rm T}Q_{1}$ with the level $\ell=13\times211$, then it is easy to verify that $S(\Delta_{1}^{\rm T})=\tilde{Q}_{1}^{\rm T}S(\Sigma)\tilde{Q}_{1}$.
   
   Similarly, for $\ell=13$, let
    \[
    Q_{2}=\frac{1}{13}\begin{bmatrix}
  -4 & 5 & 5 & 7 & -6 & 3 & 3 \\
 2 & -9 & 4 & 3 & 3 & 5 & 5 \\
 2 & 4 & -9 & 3 & 3 & 5 & 5 \\
 -4 & 5 & 5 & -6 & 7 & 3 & 3 \\
 8 & 3 & 3 & -1 & -1 & 7 & -6 \\
 8 & 3 & 3 & -1 & -1 & -6 & 7 \\
 1 & 2 & 2 & 8 & 8 & -4 & -4 \\
\end{bmatrix}.\]
   \begin{center}
       $S(\Delta_{2})=Q_2^{\rm T}S(\Sigma)Q_2=\begin{bmatrix}
       0& 1& 0& -1& 0& -1& -1\\
       -1& 0& 0& 0& 0& 1& 1\\
       0& 0&0& -1& 1&-1& 0\\
       1& 0& 1& 0& 0& -1& 1\\
       0& 0& -1& 0& 0& 1& 1\\
       1& -1& 1&1& -1& 0&-1\\
       1& -1& 0&  -1&  -1&  1&  0\\
       \end{bmatrix}.$

   \end{center}
Let $\tilde{Q}_{2}:=Q_{0}^{\rm T}Q_{2}$ with the level $\ell=7\times211$, then it is easy to verify that $S(\Delta_{2}^{\rm T})=\tilde{Q}_{2}^{\rm T}S\tilde{Q}_{2}$.
 
For $\ell=211$, let
    \[
    Q_{3}=\frac{1}{211}\begin{bmatrix}
  -13 & 29 & 138 & 121 & 23 & -96 & 9 \\
 -38 & 101 & 95 & -117 & 51 & 44 & 75 \\
 197 & 15 & 35 & -32 & 41 & -6 & -39 \\
 -5 & 141 & -93 & 79 & 90 & 28 & -29 \\
 61 & 10 & -47 & 49 & -43 & -4 & 185 \\
 -8 & -112 & 20 & 42 & 144 & 87 & 38 \\
 17 & 27 & 63 & 69 & -95 & 158 & -28 \\
\end{bmatrix}.\]
Then we have
   \begin{center}
       $S(\Delta_{3})=Q_3^{\rm T}S(\Sigma)Q_3=\begin{bmatrix}
     0& 0& 1& 0& 0&-1& 1\\
       0& 0& 1& -1& 0& 0& 1\\
       -1& -1& 0& 0& 1& 1&1\\
       0& 1& 0& 0& 1& 0& -1\\
       0& 0& -1& -1& 0& -1& 1\\
       1& 0& -1& 0&1& 0& -1\\
       -1& -1& -1& 1& -1& 1& 0\\
       \end{bmatrix}.$
   \end{center}
Let $\tilde{Q}_{3}:=Q_{0}Q_{3}$ with the level $\ell=7\times13$, then it is easy to verify that $S(\Delta_{3}^{\rm T})=\tilde{Q}_{3}^{\rm T}S\tilde{Q}_{3}$.

 Thus, we have successfully identified all cospectral mates of the $\Sigma$, namely, $\Delta_i$, $\Delta_i^{\rm T}$ ($i=1,2,3$) and $\Sigma^{\rm T}$, there are $2^{3}-1=7$ cospectral mates in total.
   \end{exam}
\section{Conclusion}\label{sec7}
In this paper, we consider a general class of oriented graphs $\mathcal{G}_{n}$, not limited to
self-converse graphs. We give a
criterion for determining whether the graph $\Sigma\in\mathcal{G}_{n}$ is $\mathrm{WDGSS}$. For non-$\mathrm{WDGSS}$ oriented graphs, we establish that the maximum number of cospectral mates of  $\Sigma\in\mathcal{G}_{n}$ is $2^{t}-1$, where $t$ is the number of prime factors of $\ell_0$.
For $t \le 3$, we have completely solved the problem of determining which graphs in $\mathcal{G}_{n}$ are $\mathrm{WDGSS}$, and of finding
out all of their generalized cospectral mates whenever they are not $\mathrm{WDGSS}$. We conclude this paper by posing the following questions for further research.
\begin{itemize}
    \item How to determine which oriented graphs in $\mathcal{G}_{n}$ are $\mathrm{WDGSS}$ for $t>3$?
    \item How to determine which oriented graphs are $\mathrm{WDGSS}$ more generally?
    \item Is it true that almost all oriented graphs are $\mathrm{WDGSS}$?
\end{itemize}

\section*{Acknowledgments}
The research of the second author is supported by National Key Research and Development Program of China 2023YFA1010203 and National Natural Science Foundation of China (Grant No. 12371357), and the third author is supported by Fundamental Research Funds for the Central Universities (Grant No. 531118010622),
National Natural Science Foundation of China (Grant No. 1240011979) and Hunan
Provincial Natural Science Foundation of China (Grant No. 2024JJ6120).

        \end{document}